\documentstyle[psfig,amscd,11pt]{amsart}
\author{Kevin P. Knudson}\thanks{Partially supported by an NSF Postdoctoral Fellowship}
\address{Department of Mathematics, Northwestern University, 
Evanston, IL  60208}
\email{knudson@@math.nwu.edu}
\date{October 7, 1997}
\title[Twisted Cohomology of $\slnft$]{Congruence Subgroups and
Twisted Cohomology of $\slnft$}

\newcommand{\slnf}{SL_n(F)}
\newcommand{\slnft}{SL_n(F[t])}

\newcommand{\gs}{\Gamma_{\sigma}}
\newcommand{\ps}{P_{\sigma}}
\newcommand{\cs}{C_{\sigma}}

\newcommand{\zz}{{\Bbb Z}}
\newcommand{\zq}{{\Bbb Q}}
\newcommand{\lra}{\longrightarrow}

\newcommand{\sps}{spectral sequence }
\newcommand{\X}{{\cal X}}

\newcommand{\T}{{\cal T}}

\newcommand{\Z}{{\cal Z}}
\newcommand{\cH}{{\cal H}}
\newcommand{\bop}{\bigoplus}

\newcommand{\ra}{\rightarrow}

\newcommand{\sslf}{SL_2(F)}
\newcommand{\sslft}{SL_2(F[t])}

\newcommand{\Hom}{{\cal H}{\em om}}
\newcommand{\lsln}{{\frak sl}_n(F)}
\newcommand{\ks}{K_{\sigma}}
\newcommand{\kv}{K_v}

\newtheorem{theorem}{Theorem}[section]
\newtheorem{prop}[theorem]{Proposition}
\newtheorem{lemma}[theorem]{Lemma}
\newtheorem{cor}[theorem]{Corollary}

\theoremstyle{remark}

\begin{document}

\maketitle

In \cite{knudson2} we showed that if $F$ is an infinite field, then
the natural inclusion $\slnf\ra \slnft$ induces an isomorphism
$$H_\bullet(\slnf,\zz) \lra H_\bullet(\slnft,\zz)$$
for all $n \ge 2$.  Here, we study the extent to which this isomorphism
holds when the trivial $\zz$ coefficients are replaced by some rational
representation of $\slnf$.  The group $\slnft$ acts on such a 
representation via the map $\slnft \stackrel{t=0}{\lra} \slnf$.

There are two approaches one might take.  The first is to use the spectral
sequence associated to the action of $\slnft$ on a certain Bruhat--Tits
building $\X$.  With this method, we obtain the following result.  Suppose
that $F$ is a field of characteristic zero.  Let $V$ be an irreducible 
representation of $\slnf$.

\noindent{ \bf Theorem.}({\em cf.} Theorem \ref{mainthm}) {\em
The group $H^1(\slnft, V)$ satisfies
$$H^1(\slnft, V) = \begin{cases}
                              H^1(\slnf, V)                               & V \ne \text{\em Ad} \\
                              H^1(\slnf, V) \oplus F^{\infty} & V=\text{\em Ad}, n=2
\\
                              H^1(\slnf, V) \oplus F  &  V=\text{\em Ad}, n\ge 3.
                             \end{cases}$$}

One expects that a similar result holds in positive characteristic as
well.  However, the author does not pretend to be an expert in 
representation theory, especially in positive characteristic; therefore,
we restrict our attention to the characteristic zero case.

The homology of linear groups with twisted coefficients has been
studied by Dwyer \cite{dwyer}, van der Kallen \cite{vdk}, and others. 
Most known results concern what happens to these homology groups as 
$n$ becomes large.  In contrast, our results cover the unstable case.

The second approach to studying this question is to use the Hochschild--
Serre spectral sequence associated to the extension
$$1 \lra K \lra \slnft \stackrel{t=0}{\lra} \slnf \lra 1$$
where $K$ denotes the subgroup of matrices which are congruent to the
identity modulo $t$.  A spectral sequence calculation shows that for any
field $F$, 
$$H^1(\slnft, V) = H^1(\slnf, V) \oplus \text{Hom}_{\slnf}(H_1(K),V)$$
so that one need only compute the group $H_1(K)$ explicitly.  This seems
to be rather difficult in general.  When $n=2$, we have a free product
decomposition
$$K=*_{s \in \sslf/B} sCs^{-1}$$
where $C$ is the upper triangular subgroup of $K$ and $B$ is the upper
triangular subgroup of $\sslf$.  Hence we have
$$H_1(K) = \bop_{s\in \sslf/B} H_1(sCs^{-1})$$
and since $C$ is an abelian group, this is easily calculated.

For $n\ge 3$, however, we have no such free product decomposition. 
In fact, the natural map $H_1(C) \ra H_1(K)$ is no longer injective 
(consider the matrix $I+E_{12}(t^2) \in C$; this gives a nontrivial
element of $H_1(C)$ which vanishes in $H_1(K)$).  The group $H_1(K)$
appears to be rather complicated in general. For any field $F$,
$H_1(K)$ surjects onto the adjoint representation ${\frak sl}_n(F)$,
and the kernel is nontrivial in general (see Section
\ref{centralseries} below).  We
make the following conjecture.

\medskip

\noindent {\bf Conjecture.} {\em  If $n\ge 3$, then $H_1(K) = {\frak sl}_n(F)$
for any finite field $F$.}

\medskip

One might prove this by finding a fundamental domain for the action of $K$
on the Bruhat--Tits building $\X$ and then utilizing the corresponding
spectral sequence.  We take the first steps toward this here.  
Unfortunately, the combinatorics involved are rather complicated.  Still,
we are able to prove the conjecture for $n=3$, $F={\Bbb F}_2, {\Bbb F}_3$.

The study of congruence subgroups has a long history, one which we will
not try to reproduce here.  The standard question one asks is whether a
given normal subgroup of $SL_n(R)$ is a congruence subgroup.  In
particular, one could ask this question for the subgroups of $K$ which
form its lower central series.  While these groups are close to being
congruence subgroups in a certain sense (see Section
\ref{centralseries}), they actually are not.

This paper consists of two parts.  The first part deals with the 
characteristic zero case and is organized as follows.
In Section \ref{sl2part1}, we
prove the above theorem in the case $n=2$ by considering the long
exact sequence associated to Nagao's amalgamated free product
decomposition of $\sslft$.  In Section \ref{spectral1}, we describe
the spectral sequence used in the proof of Theorem \ref{mainthm}.
In Section \ref{sl2part2}, we reprove the theorem in the case $n=2$
using the spectral sequence of Section \ref{spectral1}.  This helps
guide the way to the general result. In Section \ref{psigma}, we describe
some equivariant homomorphisms.  In Section \ref{h1}, we prove 
Theorem \ref{mainthm}.  

The second part deals with the congruence subgroup $K$ and its
abelianization $H_1(K)$.  In Section \ref{spectral2}, we consider the
maximal algebraic quotient of $H_1(K)$ in the characteristic zero case.
In Section \ref{centralseries} we describe a descending central series in $K$.  In
Section \ref{action} we consider the action of $K$ on the building $\X$ and
the resulting spectral sequence.  In Section \ref{reduction} we find a simpler
subcomplex of $\X$ which suffices to compute $H_1(K)$.  Finally, in
Section \ref{computation} we prove the above conjecture in certain
cases.

This paper has its origins in discussions with Dick Hain, who wondered
about the structure of $H_1(K)$ and also whether unstable homotopy
invariance holds with nontrivial coefficients.  Also, I should point out
that this work owes a great deal to C. Soul\'e's paper
\cite{soule1}. I thank Eric Friedlander for many useful discussions.
I am also grateful to the Mathematisches Forschungsinstitut Oberwolfach
and l'Institut de Recherche Math\'ematique
Avanc\'ee, Strasbourg (in particular to Jean-Louis Loday), for their hospitality during a visit in June 1996
when a preliminary version of this paper was written.  Finally, I
thank the referee who saved me from some embarrassing slips and
pointed me to Krusemeyer's paper \cite{kruse}.

\medskip

\noindent {\em Notation.}  In Sections 1 through 6, $F$ is assumed to be
of characteristic zero; thereafter, $F$ is allowed to be any field
unless otherwise specified.  
If $G$ is a group acting on a set X, we denote the
invariants of the $G$--action by $X^G$.  The symbol $P$ typically
denotes a parabolic subgroup of $\slnf$ which contains the upper
triangular subgroup $B$.  We denote by $\Phi^+$ the set of positive roots
of $\slnf$ determined by $B$; this set consists of roots $\{ \alpha_{ij}:
i<j \}$.  We denote by $E_{ij}(a)$ the matrix having $i,j$
entry $a$ and zeros elsewhere.  For a polynomial $p(t)$, $p^{(k)}$ denotes
the coefficient of $t^k$.  If $R$ is a ring, we denote by $R^\times$
the group of units of $R$.

\section{The $SL_2$ Case, Part I}\label{sl2part1}

We single out the case $n=2$ because we may use Nagao's Theorem
\cite{nagao}
\begin{equation}\label{amal}
\sslft = \sslf *_B B_t
\end{equation}
to study $H^1(\sslft,V)$ (here, $B_t$ denotes the group of upper triangular matrices over $F[t]$).  This amalgamated free product decomposition
yields a long exact sequence
$$\cdots \ra H^k(\sslft,V) \ra H^k(\sslf,V)\oplus H^k(B_t,V) \ra H^k(B,V)
\ra \cdots$$
for computing the cohomology of $\sslft$.

\begin{prop}
We have a short exact sequence
\begin{equation}\label{h1sl2ses}
0\ra H^1(\sslft,V) \ra H^1(\sslf,V) \oplus H^1(B_t,V) \ra
H^1(B,V) \ra 0.
\end{equation}
\end{prop}

\begin{pf}  Since the map $B_t \stackrel{t=0}{\lra} B$ is split by the
natural inclusion $B\lra B_t$, the induced map $H^k(B_t,V) \lra H^k(B,V)$
is surjective.
\end{pf}

We now study the relationship between $H^1(B_t,V)$ and $H^1(B,V)$.  
Consider the extension
\begin{equation}\label{bext}
1 \lra C \lra B_t \stackrel{t=0}{\lra} B \lra 1
\end{equation}
where
$$C=\biggl\{ \left(\begin{array}{cc}
                                 1  &  tp(t) \\
                                 0  &   1
                               \end{array} \right): p(t) \in F[t] \biggr\}.$$
The Hochschild--Serre spectral sequence associated to this
satisfies
\begin{equation}\label{ss1}
E_2^{p,q} = H^p(B,H^q(C,V)) \Longrightarrow H^{p+q}(B_t,V).
\end{equation}
Since the extension (\ref{bext}) is split, the map $d_2^{0,1}: E_2^{0,1}
\ra E_2^{2,0}$ vanishes.  It follows that we have a short exact sequence
\begin{equation}\label{bth1}
0\ra H^1(B,H^0(C,V)) \ra H^1(B_t,V) \ra H^0(B,H^1(C,V))\ra 0.
\end{equation}
Since $C$ acts trivially on $V$, $H^0(C,V) = V$.  Hence, the first term
in (\ref{bth1}) is simply $H^1(B,V)$.  Observe that
$H^0(B,H^1(C,V))=\text{Hom}_B(H_1(C),V)$.

\begin{prop}\label{laterprop}
The group $\text{\em Hom}_B(H_1(C),V)$ satisfies
$$\text{\em Hom}_B(H_1(C),V)= \begin{cases}
                                   \text{\em Hom}_F(tF[t],F) & V=\text{\em Ad} \\
                                                    0           & V \ne \text{\em Ad}.
                               \end{cases}$$
\end{prop}

\noindent {\em Proof.}  
We may compute this group as follows.  We have a split extension
$$1\lra U\lra B\lra T\lra 1$$
where $U$ is the subgroup of upper triangular unipotent matrices and 
$T$ is the diagonal subgroup.  The Hochschild--Serre \sps
implies that
$$H^0(B,H^1(C,V)) = H^0(T,H^0(U,H^1(C,V)));$$
that is,
$$\text{Hom}_B(H_1(C),V) = H^0(T,\text{Hom}_U(H_1(C),V)).$$
Observe that $H_1(C) = C \cong tF[t]$ and $U$ acts trivially on $C$.  It
follows that
$$\text{Hom}_U(H_1(C),V) = \text{Hom}(C,V^U).$$

Now, $V=\text{Sym}^n S$, where $S$ is the standard representation of 
$\sslf$ and $n$ is some nonnegative integer (see {\em e.g.} \cite[Ch. 11]{fulton}).  It follows that
$$V^U = V_n$$ where $V_n$ denotes the highest weight space of $V$.  
Observe that $T$ acts on $C$ with weight 2 and on $V_n$ with weight
$n$.  Note also that
$$H^0(T,\text{Hom}(C,V_n))  = \text{Hom}_T(C,V_n).$$
We need the following result.

\begin{lemma}\label{lin2}  Suppose $f:C \lra V^U$ is a $T$-equivariant
group homomorphism.  Then $f$ is $F$-linear.
\end{lemma}

\noindent {\em Proof.}  We have $V=\text{Sym}^nS$ and $V^U = V_n$.  If
$a \in \zz$, denote by $t_a$ the element of $T$ with diagonal entries
$a, \frac{1}{a}$.  If $w \in C$, we have $t_a.w = a^2w$.  Thus
$$a^nf(w) = t_a.f(w) = f(t_a.w) = f(a^2w) = a^2f(w)$$
(the last equality follows since $f$ is $\zz$-linear).  It follows
that $f=0$ unless $n=2$.

Consider the case $n=2$ ({\em i.e.,} $V=\text{Ad}$).  Let $\alpha \in
F$ and consider the element $t_\alpha = \text{diag}(\alpha,
\frac{1}{\alpha})$.  Since $f$ is $T$-equivariant, we have
$$\alpha^2f(w) = t_\alpha.f(w) = f(t_\alpha.w) = f(\alpha^2w)$$
for any $w \in C$.  Note that since $F$ has characteristic zero and
$f$ is $\zz$-linear, $f$ is in fact $\zq$-linear.  Note also that for
any $\alpha \in F$, we have $\alpha = -\frac{1}{2} -
\frac{1}{2}\alpha^2 + \frac{1}{2}(1+\alpha)^2$.  Thus,
\begin{eqnarray*}
f(\alpha w) & = & f((-\frac{1}{2} - \frac{1}{2}\alpha^2 +
\frac{1}{2}(1+\alpha)^2)w) \\
            & = & -\frac{1}{2}f(w) - \frac{1}{2}\alpha^2f(w) +
            \frac{1}{2}(1+\alpha)^2f(w) \\
            & = & \alpha f(w). 
\end{eqnarray*} \hfill \qed
    
\medskip

Returning to the proof of Proposition \ref{laterprop}, the proof of
Lemma \ref{lin2} shows that $\text{Hom}_T(C,V_n)=0$ unless $n=2$.  It
follows that 
$$\text{Hom}_B(C,V) = \begin{cases}
                                                0                  &  V \ne \text{Ad}  \\
                                      \text{Hom}_F(tF[t], F)       &  V = \text{Ad}.  
                                     \end{cases}$$ \hfill \qed 

\medskip

\begin{cor}
The group $H^1(B_t,V)$ satisfies
$$H^1(B_t,V) = \begin{cases}
                            H^1(B,V)                        &  V \ne \text{\em Ad}  \\
                            H^1(B,V) \oplus \text{\em Hom}_F(tF[t],F)  &  V= \text{\em Ad}.
                         \end{cases}$$
\end{cor}

\begin{pf}  This follows from the short exact sequence (\ref{bth1}) and
the fact that $H^1(B,V)$ is a direct summand of $H^1(B_t,V)$.
\end{pf}

\begin{theorem}\label{sl2h1}
The group $H^1(\sslft,V)$ satisfies
$$H^1(\sslft, V) = \begin{cases}
                                 H^1(\sslf, V)             & V\ne \text{\em Ad}  \\
                                 H^1(\sslf,V)\oplus \text{\em Hom}_F(tF[t],F) & V=\text{\em Ad}. 
\end{cases}$$
\end{theorem}

\begin{pf}  Consider the short exact sequence (\ref{h1sl2ses}).  By the
corollary, the kernel of the map
$$H^1(\sslf,V)\oplus H^1(B_t,V) \lra H^1(B,V)$$
is isomorphic to
$$\begin{array}{cl}
  H^1(\sslf,V)               &    V \ne \text{Ad}  \\
H^1(\sslf,V) \oplus \text{Hom}_F(tF[t],F) &  V = \text{Ad}.
\end{array}$$
The direct sum decomposition follows because $H^1(\sslf,V)$ is a direct
summand of $H^1(\sslft,V)$.
\end{pf}

We have seen that $H^1(\sslft,\text{Ad})$ differs from $H^1(\sslf,\text{Ad})$ by the infinite dimensional $F$-vector space
$\text{Hom}_B(H_1(C),\text{Ad})$.  We describe an explicit basis for
this space.  For each $k\ge 1$, denote by $\varphi_k$ the map
$C\ra {\frak sl}_2(F)$ defined by
$$I + E_{12}(a_{12}^{(1)}t + a_{12}^{(2)}t^2 + \cdots + a_{12}^{(m)}t^m) \mapsto
E_{12}(a_{12}^{(k)}).$$
Then the set $\{\varphi_k\}$ is a basis of $\text{Hom}_B(H_1(C),\text{Ad})$.

\section{A Spectral Sequence}\label{spectral1}

In general, we shall use the action of $\slnft$ on a suitable simplicial
complex to compute $H^1(\slnft,V)$.

Denote by $\X$ the Bruhat--Tits building associated to the vector space
$F(t)^n$.  Recall that the vertices of $\X$ are equivalence
classes of ${\cal O}$-lattices in $F(t)^n$ (here, ${\cal O}$ consists of
the set of $a/b$ with $\deg b \ge \deg a$), where two lattices $L$ and
$L'$ are equivalent if there is an $x\in F^\times$ with $L'=xL$.
  A collection of vertices
$\Lambda_0,\Lambda_1,\dots ,\Lambda_m$ forms an $m$-simplex if
there are representatives $L_i$ of the $\Lambda_i$ with
$$t^{-1}L_0 \subset L_m \subset L_{m-1} \subset \cdots \subset L_0.$$
It is possible to put a metric on $\X$ so that each edge in $\X$ has
length one.  When we speak of the distance between vertices we implicitly
use this metric.
For a more complete description of $\X$, see, for example, \cite{knudson2}.

 The
group $\slnft$ acts on $\X$ with fundamental domain an infinite wedge
$\T$, which is the subcomplex of $\X$ spanned by the vertices
$$[e_1t^{r_1},e_2t^{r_2},\dots , e_{n-1}t^{r_{n-1}}, e_n], \quad r_1 \ge r_2 \ge
\cdots r_{n-1} \ge 0$$ where $e_1,e_2,\dots ,e_n$ denotes the
standard basis of $F(t)^n$ 
(this is due to Serre \cite{serre} for $n=2$ and Soul\'e \cite{soule1} for
$n\ge 3$).  See Figure 1 for the case $n=3$.

\begin{figure}
\centerline{\psfig{figure=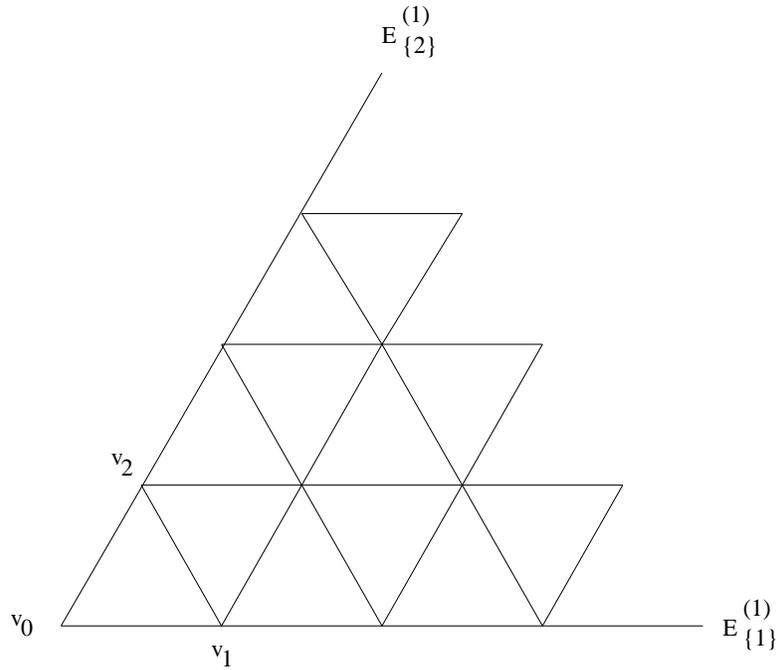,angle=270,height=3.5in}}
\caption{The Fundamental Domain $\T$ for $n=3$}
\end{figure}

Denote by $v_0$ the vertex $[e_1, \dots,
e_n]$ and by $v_i$ the vertex $$[e_1t,e_2t, \dots ,e_it,e_{i+1}, \dots, e_n],\quad
i=1,2,\dots,n-1.$$  For a $k$ element subset $I=\{i_1, \dots , i_k\}$ of 
$\{1,2, \dots , n-1\}$, define $E_I^{(k)}$ to be the subcomplex of $\T$ which
is the union of all rays with origin $v_0$ passing through the $(k-1)$--simplex
$\langle v_{i_1}, \dots , v_{i_k}\rangle$.  There are $\binom {n-1}k$                                                                                                                  such $E_I^{(k)}$.  Observe that if $I=\{1,2, \dots , n-1\}$, then
$E_I^{(n-1)} = \T$.  When we write $E^{(l)}_J$, the superscript $l$ denotes the
cardinality of the set $J$.

The structure of the various simplex stabilizers was described in
\cite{knudson2}.  If $(x_{ij}(t)) \in \slnft$ stabilizes the vertex
$[e_1t^{r_1},\dots e_{n-1}t^{r_{n-1}},e_n]$, then we have
$$\deg x_{ij}(t) \le r_i - r_j$$ (set $r_n=0$).  Note that since $r_1\ge r_2 \ge
\cdots \ge r_{n-1}$, some of the $x_{ij}(t)$ with $i>j$ may be $0$.
 Denote the stabilizer of $\sigma$ by
$\gs$.  The group $\gs$ is the intersection of the stabilizers $\Gamma_v$ where $v$ ranges over the vertices of $\sigma$.  In this case $\deg
x_{ij}(t) \le \displaystyle \min_{v \in \sigma}\{r_i^{(v)} - r_j^{(v)}\}$.
Observe that in any case, the group $\gs$ has a block form where the
blocks along the diagonal are matrices with entries in $F$, blocks below 
are zero, and blocks above contain polynomials of bounded degree.
In the case $n=3$, we have the following block forms:
$$\begin{array}{rcl}
v_0                                &  :  &  \Gamma_{v_0} = SL_3(F)                   \\
\sigma \in E^{(1)}_{\{1\}} & : & {\gs=\left(\begin{array}{c|cc}
                                                                    *  &  *  &   *   \\ \hline
                                                                    0  &  *  &   *    \\
                                                                    0  &  *  &   *
                                                           \end{array}\right)  }                 
\end{array}$$

$$\begin{array}{rcl}
\sigma \in E^{(1)}_{\{2\}} & : & {\gs=\left(\begin{array}{cc|c}
                                                                    *  &  *  &   *   \\    
                                                                     *  &  *  &   *    \\ \hline
                                                                    0  &  0  &   *
                                                           \end{array}\right)  }                 \\
\sigma \in \T - (E^{(1)}_{\{1\}} \cup E^{(1)}_{\{2\}}) & : & 
                                                   {\gs=\left(\begin{array}{ccc}
                                                                    *  &  *  &   *   \\ 
                                                                    0  &  *  &   *    \\
                                                                    0  &  0  &   *
                                                           \end{array}\right)  } 
\end{array}$$

We have a short exact sequence
$$1 \lra C_{\sigma} \lra \gs \stackrel{t=0}{\lra} \ps \lra 1$$
where $\ps$ is a parabolic subgroup of $\slnf$.   From the above 
description of $\gs$, we see that the group $\cs$ has a block form where
blocks along the diagonal are identity matrices, blocks below are zero,
and blocks above contain polynomials of bounded degree which are
divisible by $t$.

Filter
the complex $\T$ by setting $W^{(0)} = v_0$ and
$$W^{(l)} = \bigcup_{|I|=l} E_I^{(l)}, \; \; 1\le l \le n-1$$
  Observe that if $\sigma$ and $\tau$ are 
simplices in the same component of $W^{(i)} - W^{(i-1)}$, then
$\ps = P_{\tau}$ since on such a component the relationships among
the $r_i$ defining the vertices do not vary from vertex to vertex
({\em i.e.,} if $r_i>r_{i+1}$ for one vertex in the component, then the
same holds for every vertex in the component; since these relations
determine which entries below the diagonal are zero, we see that the
stabilizers of these vertices all have the same block form and hence so
does the stabilizer of any simplex in the component).

The action of $\slnft$ on $\X$ gives rise to a spectral sequence converging to $H^{\bullet}
(\slnft, V)$ with $E_1$--term
\begin{equation}\label{ss2}
E_1^{p,q} = \prod_{\dim\sigma = p} H^q(\gs, V)
\end{equation}
where $\gs$ denotes the stabilizer of the simplex $\sigma \in \T$.  We
shall compute the terms $E_2^{p,0}$ and $E_2^{0,1}$ and thereby
obtain the group $H^1(\slnft,V)$.

\begin{prop}
The bottom row of the spectral sequence {\em (\ref{ss2})}
satisfies $$E_2^{p,0} = 0, \quad p> 0.$$
\end{prop}

\begin{pf}  We have a coefficient system ${\cal H}^0$ on $\T$ given by
$${\cal H}^0(\sigma) = H^0(\gs, V).$$  Thus, the bottom row of (\ref{ss2})
is simply the cochain complex  $C^{\bullet}(\T, {\cal H}^0)$.  The Hochschild--Serre \sps associated to the extension $\gs \ra \ps$ implies
that
$$H^0(\gs , V) = H^0(\ps , V).$$
It follows that on each component of $W^{(i)}-W^{(i-1)}$, the
system ${\cal H}^0$ is constant (see the remarks above).  An easy modification of Lemma 5 of 
\cite{soule2} shows that the inclusion $v_0 \lra \T$ induces an
isomorphism
$$H^{\bullet}(\T, {\cal H}^0) \stackrel{\cong}{\lra} H^{\bullet}(v_0, {\cal H}^0)$$
(see \cite{knudson2} for the corresponding statement for homology).
It follows that
$$E_2^{p,0} = \begin{cases}
                           H^0(\slnf, V)   &   p=0  \\
                                  0              &   p>0.
                        \end{cases}$$
\end{pf}

It remains to compute the term 
$$E_2^{0,1} = \ker\{ d_1:E_1^{0,1}\lra E_1^{1,1}\}$$
in the spectral sequence (\ref{ss2}).  That is, we must compute the
group $H^0(\T, {\cal H}^1)$, where ${\cal H}^1$ is the system which 
assigns the group $H^1(\gs, V)$ to $\sigma$.

\section{The $SL_2$ Case, Part II}\label{sl2part2}

Before computing the group $H^0(\T, {\cal H}^1)$ in general, we first
reconsider the case $n=2$.  This will help guide the way to the desired
result.  To this end, consider the fundamental domain $\T$.  In this
case, $\T$ is an infinite path in the tree $\X$.  Label the vertices
$v_0,v_1, \dots$ and the edges $e_0,e_1,\dots$.  Then we have 
\begin{eqnarray*}
\Gamma_{v_0}   &  =  &  \sslf  \\
\Gamma_{v_i}    &  =  &  {\Gamma_{e_i} = \biggl\{ \left(\begin{array}{cc}
                                                                                 a & p(t) \\
                                                                                 0 & 1/a
                                                                                 \end{array}\right):
                                                \deg p(t) \le i \biggr\}, \; i \ge 1.}
\end{eqnarray*}
For each $i\ge 1$ we have a split extension
$$1\lra C_i \lra \Gamma_{v_i} \stackrel{t=0}{\lra} B \lra 1.$$
Arguing as in the proof of Proposition \ref{laterprop} and again using
Lemma \ref{lin2} we see that
$$H^1(\Gamma_{v_i}, V) = \begin{cases}
                                              H^1(B,V)                   &  V \ne \text{Ad} \\
                                              H^1(B,V) \oplus \text{Hom}_B(C_i,V) & V=
                                                                                   \text{Ad}.
                                          \end{cases}$$

\begin{theorem}
The group $H^0(\T, {\cal H}^1)$ satisfies
$$H^0(\T,{\cal H}^1) = \begin{cases}
                                      H^1(\sslf, V)                & V\ne \text{\em Ad} \\
                                      H^1(\sslf, V) \oplus \text{\em Hom}_B(C,V) & 
                                                                                V = \text{\em Ad}.
                                    \end{cases}$$
\end{theorem}

\begin{pf}  If $V$ is not the adjoint representation, we use the long exact
sequence of the pair $(\T,v_0)$.  Since $H^1(\Gamma_{v_i},V) = 
H^1(B,V)$ for each $i\ge 1$, we see that $H^\bullet(\T,v_0;\cH^1) = 0$.  It follows that
$$H^0(\T, {\cal H}^1) = H^1(\sslf, V).$$  In the case $V=\text{Ad}$, we
have
$$C^{\bullet}(\T,{\cal H}^1) = C^{\bullet}(\T,{\cal H}^1(P,V)) \oplus
                                                  C^{\bullet}(\T,\Hom(C,V))$$
where ${\cal H}^1(P,V)$ is the system
$$\sigma \mapsto H^1(\ps, V)$$
and $\Hom(C,V)$ is the system
$$\sigma \mapsto \text{Hom}_{\ps}(C_{\sigma},V).$$

Again, we see that $H^0(\T,{\cal H}^1(P,V)) = H^1(\sslf, V)$.  It remains
to compute $H^0(\T, \Hom(C,V))$.  For each $i$, we have
\begin{eqnarray*}
\text{Hom}_B(C_i,V)   &  =   &  \text{Hom}_T(\{a_{12}^{(1)}t + \cdots
                                                    + a_{12}^{(i)}t^i\}, V) \\
                                   &  =   & F\{\varphi_j\}_{j=1}^i
\end{eqnarray*}
where $\varphi_j(a_{12}^{(1)}t+\cdots +a_{12}^{(i)}t^i) = 
E_{12}(a_{12}^{(j)})$ (this follows from Lemma \ref{lin2}).  Note that $C_{v_0}=C_{e_0}=1$; {\em i.e.}
$\text{Hom}(C_{v_0},V)=0$.  It follows that $$H^0(\T,\Hom(C,V)) =
H^0(\T',\Hom(C,V)),$$ where $\T'$ is the complex obtained from
deleting $v_0$ and $e_0$ (but not $v_1$) from $\T$.

Consider the map $d:C^0(\T',\Hom(C,V)) \ra C^1(\T', \Hom(C,V))$.
This map is given by
$$(x_1,x_2, x_3, \dots ) \mapsto (x_2|_{e_1}-x_1,x_3|_{e_2}-x_2, \dots)
$$  where $x_k|_{e_{k-1}}$ denotes the image of $x_k$ under the
restriction map $$\text{Hom}_B(C_{v_k},V) \ra \text{Hom}_B(C_{e_{k-1}},
V).$$  Note that $d$ is $F$-linear since it is an alternating sum of
restriction maps, each of which is clearly $F$-linear.
The kernel of this map is spanned by the following elements:
$$(\varphi_1, \varphi_1, \varphi_1, \dots ),
(0, \varphi_2, \varphi_2,\varphi_2, \dots ),$$
$$(0,0, \varphi_3, \varphi_3, \varphi_3, \dots),
\dots $$
If we identify the $k$th element of this list with the map $\varphi_k:C \ra V$,
we see that
$$H^0(\T',\Hom(C,V)) = \text{Hom}_B(C,V).$$  It follows that
$$H^1(\sslft, \text{Ad}) = H^1(\sslf, \text{Ad}) \oplus \text{Hom}_B(C,
\text{Ad})$$ which is the desired result.
\end{pf}

\section{$\ps$--equivariant Homomorphisms}\label{psigma}

\newcommand{\vak}{v_\alpha^{(k)}}

To compute the group $H^0(\T,{\cal H}^1)$ in general, we will need to 
consider the groups $\text{Hom}_{\ps}(H_1(\cs),V)$, where the groups
$H_1(\cs)$ are more complicated than in the case $n=2$.  We begin by
describing the structure of the $H_1(\cs)$.

Recall that the elements $(x_{ij}(t))$ of the subgroup $\cs$ have a block form where
the blocks on the diagonal are identity matrices, blocks below are zero,
and blocks above contain polynomials of bounded degree which are divisible by $t$. 
 
Using elementary row operations, we may write for any $(x_{ij}(t)) \in \cs$, 

\begin{equation}\label{prod1}
(x_{ij}(t)) = \prod_{j=0}^{n-2} \prod_{i=1}^{n-j-1} (I + E_{i,n-j}(x_{i,n-j}(t))).
\end{equation}
We also have

\begin{equation}\label{prod2}
I + E_{ij}(x_{ij}(t)) = \prod_{l=1}^m (I + E_{ij}(x_{ij}^{(l)} t^l)).
\end{equation}

\begin{lemma}\label{comm}
If there exist distinct integers $i,j,k$ and positive integers $l,m$ with
$l+m = r$ such that $I+E_{ik}(t^l), I+E_{kj}(t^m) \in \cs$, then the element
$I+E_{ij}(t^r)$ is zero in $H_1(\cs)$.
\end{lemma}

\begin{pf}  This follows from standard relations among commutators of
elementary matrices.
\end{pf}

We are now in a position to describe explicitly the structure of the groups
$H_1(\cs)$.

\begin{prop}  The group $H_1(\cs)$ admits a decomposition
$$H_1(\cs) = \bop_{\alpha \in \Phi^+} S_\alpha$$
where $S_\alpha$ is a weight space (possibly 0) for the positive root
$\alpha$.  Moreover, each $S_{\alpha_{ij}}$ is graded by powers of $t$:
$$S_{\alpha_{ij}} = t^{l_1(\alpha_{ij})}W_{\alpha_{ij}} \oplus t^{l_2(\alpha_{ij})}W_{\alpha_{ij}} \oplus \cdots \oplus t^{l_p(\alpha_{ij})}W_{\alpha_{ij}}$$
where $0\le l_1(\alpha_{ij}) < l_2(\alpha_{ij}) < \cdots < l_p(\alpha_{ij})
\le \displaystyle \min_{v \in \sigma} \{r_i^{(v)} - r_j^{(v)}\}$ and 
$W_{\alpha_{ij}}$ is a one-dimensional root space for $\alpha_{ij}$.
If $S_{\alpha_{ij}} \ne 0$, then $l_1(\alpha_{ij}) = 1$.
\end{prop}

\begin{pf}
Observe that $H_1(\cs)$ is a $\ps$--module and as such admits a 
decomposition into weight spaces for the action of the diagonal subgroup
$T$.  In light of (\ref{prod1}) and (\ref{prod2}), we see that each element
of $H_1(\cs)$ may be written as a product
$$\prod_{i<j} I+E_{ij}(x_{ij}(t))$$
which implies the existence of the decomposition
$$H_1(\cs) = \bop_{\alpha\in \Phi^+} S_\alpha$$
(the only weights which can occur are positive roots since $\cs$ is
upper triangular and $T$ acts on $I+E_{ij}(x)$ with weight $\alpha_{ij}$).
The space $S_{\alpha_{ij}}$ is spanned
by  all elements of the form
$$I+E_{ij}(t^r)$$
where $1\le r \le \displaystyle \min_{v \in \sigma}\{r_i^{(v)} - r_j^{(v)}\}$. If $\displaystyle \min_{v\in \sigma}\{r_i^{(v)} - r_j^{(v)}\} =
0$, then $S_{\alpha_{ij}} = 0$. By Lemma \ref{comm}, if there is a $k$ distinct from $i$ and
$j$ with $I+E_{ik}(t^l),I+E_{kj}(t^m) \in \cs$ ($l+m = r$), then $I+E_{ij}(t^r) = 0$ in $H_1(\cs)$.  Notice that the conditions of Lemma
\ref{comm} cannot be met for $r=1$ (since $l,m\ge 1$) so that if $S_{\alpha_{ij}} \ne 0$ (which can only occur if $\displaystyle \min_{v\in \sigma}\{r_i^{(v)} - r_j^{(v)}\}\ge 1$), then $I+E_{ij}(t) \in S_{\alpha_{ij}}$.  It follows
that
$$S_{\alpha_{ij}} = tW_{\alpha_{ij}} \oplus t^{l_2(\alpha_{ij})}W_{\alpha_{ij}} \oplus \cdots \oplus t^{l_p(\alpha_{ij})}W_{\alpha_{ij}}$$
where $1 < l_2 < \cdots < l_p \le \displaystyle \min_{v \in \sigma} \{r_i^{(v)} - r_j^{(v)}\}$.
\end{pf}

Note that it is possible for powers of $t$ to get skipped.  For example, let
$n=3$ and consider the following group:
$$\cs = \left(\begin{array}{ccc}
                       1 & \deg \le 1 & \deg \le 3 \\
                       0 &        1        & \deg \le 1 \\
                       0 &         0        &    1
                        \end{array} \right).$$
This is the $\cs$ of some edge in $\T$.  In $H_1(\cs)$, the element
$$\left(\begin{array}{ccc}
              1  &  0  & t^3 \\
              0  &  1  &  0  \\
              0  &  0  &  1
             \end{array} \right)$$
is nonzero, but the element
$$\left(\begin{array}{ccc}
              1  &  0  & t^2 \\
              0  &  1  &  0  \\
              0  &  0  &  1
             \end{array} \right)$$
is trivial.  (In this case, $H_1(\cs) = tW_{\alpha_{12}} \oplus tW_{\alpha_{23}} \oplus tW_{\alpha_{13}} \oplus t^3W_{\alpha_{13}}$.)
This problem arises because for each $i,j$ we are only allowed polynomials of degree at most $\displaystyle \min_{v\in\sigma}\{r_i^{(v)} 
- r_j^{(v)}\}$.  This may allow certain elements to survive in $H_1$.

However, if $v$ is a vertex, this phenomenon does not occur for the following reason.  For each
$i,j$,  we have $\deg x_{ij}(t) \le r_i-r_j$. Consider the elements
$I+E_{ij}(t^r)$, $1\le r \le r_i-r_j$, in $C_v$.  Suppose there exists
a $k$ such that $C_v$ contains nonzero elements of the form
$I+E_{ik}(p(t))$, $I+E_{kj}(q(t))$.  Note that since $C_v$ is upper 
triangular, such a $k$ satisfies $i<k<j$.  Then $C_v$ contains all
the $I+E_{ik}(t^l)$, $I+E_{kj}(t^m)$, $1\le l \le r_i-r_k$, $1\le m \le r_k
-r_j$.  Since $r_i\ge r_k \ge r_j$, it follows that it is always possible
to satisfy the conditions of Lemma \ref{comm} for each $1< r \le r_i-r_j$;
{\em i.e.,} $S_{\alpha_{ij}} = tW_{\alpha_{ij}}$.  Note that such a pair
$i,j$ must satisfy $j>i+1$ since each $C_v$ is upper triangular.

We now determine the structure of the groups $\text{Hom}_{\ps}(H_1(\cs),V)$.

\begin{lemma}\label{lin3}
Suppose $n\ge 3$ and $f:H_1(\cs) \lra V$ is a $\ps$-equivariant group
homomorphism. Then $f(W_{\alpha_{ij}}) \subseteq V_{\alpha_{ij}}$.
\end{lemma}

\begin{pf}  Write $V=\bop_{\lambda} V_\lambda$ and choose a basis
$\{v_\lambda\}_\lambda$ of $V$.  Recall that each $\lambda$ is a
linear combination $\lambda = \sum_{p=1}^n m_p(\lambda) L_i$ where
$L_1,\dots , L_n$ are the weights of the standard representation.
Recall also that $\sum_{i=1}^n L_i = 0$.  Fix $i<j$.  For $k\ne i$
denote by $t_k$ the element $$t_k = \text{diag}(1,\dots ,1, b, 1,
\dots , 1, 1/b, 1, \dots , 1)$$ where $b$ is an integer greater than
$1$, $b$ appears in the $i$th place and $1/b$ appears in the $k$th
place.  Let $w \in W_{\alpha_{ij}}$ and write $f(w)= \sum_\lambda
a_\lambda v_\lambda$.  Note that if $k\ne j$, $t_k.w = bw$ and
$t_j.w = b^2w$.  Now if $k\ne j$,
\begin{eqnarray*}
t_k.\sum_\lambda a_\lambda v_\lambda & = & t_k.f(w) \\
       &  =  & f(t_k.w) \\
       &  =  & f(bw) \\
       &  =  & bf(w) \qquad (\text{$f$ additive}) \\
       &  =  & \sum_\lambda ba_\lambda v_\lambda
\end{eqnarray*}
and similarly,
$$t_j.\sum_\lambda a_\lambda v_\lambda = \sum_\lambda b^2a_\lambda
v_\lambda.$$
Now, if $\lambda = \sum_{p=1}^n m_p(\lambda)L_p$, then
$$t_k.v_\lambda = b^{m_i(\lambda) - m_k(\lambda)}v_\lambda$$
for all $k\ne i$.  Thus,
$$(b-b^{m_i(\lambda) - m_k(\lambda)})a_\lambda = 0, \qquad k \ne j$$
and 
$$(b^2 - b^{m_i(\lambda) - m_j(\lambda)})a_\lambda = 0.$$

So, if some $m_i(\lambda) - m_k(\lambda) \ne 1$ or $m_i(\lambda) -
m_j(\lambda) \ne 2$, we see that $a_\lambda = 0$ (since $F$ has
characteristic zero).

Suppose that for $k\ne j$, $m_i(\lambda) - m_k(\lambda) = 1$ and
$m_i(\lambda) - m_j(\lambda) = 2$.  Then
\begin{eqnarray*}
\lambda & = & (m_i(\lambda) - 1)L_1 +\cdots + m_i(\lambda)L_i + \cdots
+ (m_i(\lambda) - 2)L_j + \cdots \\
        &   & \cdots + (m_i(\lambda) - 1)L_n \\
  &  =  &  m_i(\lambda)\sum_{p=1}^nL_p - \sum_{p\ne i}L_p - L_j \\
  &  =  &  m_i(\lambda)\cdot 0 + L_i - \sum_{p=1}^n L_p - L_j \\
  &  =  &  0 + L_i - 0 - L_j \\
  &  =  &  \alpha_{ij}.
\end{eqnarray*}

Thus, $a_\lambda = 0$ except for $\lambda = \alpha_{ij}$; {\em i.e.,}
$f(w)= a_{\alpha_{ij}}v_{\alpha_{ij}}$.  
\end{pf}

\begin{cor}\label{lin3cor}  If $n\ge 3$ and $f:H_1(\cs)\lra V$ is a 
$\ps$-equivariant group homomorphism, then $f$ is $F$-linear.
\end{cor}

\noindent {\em Proof.}  Let $w \in W_{\alpha_{ij}}$ and let $\alpha
\in F$.  Let $$t_\alpha = \text{diag}(1,\dots , 1, \alpha, 1, \dots ,
1, 1/\alpha, 1, \dots , 1)$$ where $\alpha$ appears in the $i$th
position and $1/\alpha$ appears in the $k$th position, where $k\ne j$.
Then $t_\alpha.w = \alpha w$.  Hence,
$$t.f(w) = f(t.w) = f(\alpha w).$$
Since $f(W_{\alpha_{ij}})\subseteq V_{\alpha_{ij}}$, $t.f(w) = \alpha
f(w)$.  Thus, $f(\alpha w) = \alpha f(w)$ for all $\alpha \in F$.
Since each element of $H_1(\cs)$ is an $F$-linear combination of
elements of the various $W_{\alpha_{ij}}$, we see that $f$ is $F$-linear.
\hfill $\qed$

The preceding result allows us to assume that the homomorphisms under
consideration are $F$-linear.

\begin{prop}\label{vanishing} If $V$ is not the adjoint representation, then
$$\text{\em Hom}_{\ps}(H_1(\cs),V)=0.$$
\end{prop}

\begin{pf} Suppose $f$ is a nonzero element in $\text{Hom}_{\ps}
(H_1(\cs),V)$.  Since $V$ is irreducible, the image of $f$ contains a
highest weight vector whose weight $\gamma$ must be a weight of
$H_1(\cs)$, hence a root.  On the other hand, $\gamma$ must be a
highest weight of $V$, so we must have $V=\text{Ad}$.
\end{pf}

This argument also shows that if $V=\text{Ad}$, then the restriction
map $\text{Hom}_{\ps}(H_1(\cs),V) \ra \text{Hom}_F(H_1(\cs)_{\alpha_{1n}},V_{\alpha_{1n}})$ is injective.
Observe that  each map
$H_1(\cs)_{\alpha_{1n}} \ra V_{\alpha_{1n}}$ is a linear combination
of maps of the form
$$I + E_{1n}(a_{1n}^{(1)}t + a_{1n}^{(2)}t^2 + \cdots + a_{1n}^{(l)}t^l) \mapsto
E_{1n}(a_{1n}^{(k)}).$$  It follows that any $f$ in $\text{Hom}_{\ps}(H_1(\cs),V)$ is a linear combination of maps
$\varphi_k:H_1(\cs) \ra V$ defined by
$$\varphi_k:I+tX_1 +\cdots + t^lX_l \mapsto X_k.$$

Let us consider what the proposition says about $\text{Hom}_B
(H_1(C),\text{Ad})$.  If $v$ is a vertex in $\T - W^{(n-2)}$
({\em i.e.,} $v$ is in the interior of $\T$), then
$$H_1(C_v) = \bop_i W^{n_{\alpha_{i,i+1}}}_{\alpha_{i,i+1}} \oplus \bop_{j>i+1}tW_{\alpha_{ij}}$$
where $n_{\alpha_{i,i+1}}\ge 1$.  In this case, $P_v = B$ and $\text{Hom}_{P_v}(H_1(C_v),\text{Ad})$ injects into
$\text{Hom}_F(H_1(C_v)_{\alpha_{1n}}, \text{Ad}_{\alpha_{1n}})
= F$ (we assume that $n\ge 3$).  Hence, up to scalars any nonzero map
$f:H_1(C_v) \ra \text{Ad}$ is of the form
$$f(I+tX_1+t^2X_2+\cdots +t^kX_k) = X_1.$$
Note that when $n=2$ the situation is much different since $W_{\alpha_{12}}^{n_{\alpha_{12}}}$ is such that $n_{\alpha_{12}}$
can be greater than 1.

Thus we arrive at the following description of  $\text{Hom}_B(H_1(C),\text{Ad})$.

\begin{cor} The group $\text{\em Hom}_B(H_1(C),\text{\em Ad})$ satisfies
$$\text{\em Hom}_B(H_1(C),\text{\em Ad}) = \begin{cases}
                                                            F\{\varphi_k\}_{k=1}^\infty & n=2\\
                                                            F\{\varphi_1\} & n\ge 3
                                                         \end{cases}$$
where $\varphi_k$ is the map $H_1(C) \lra \text{\em Ad}$ defined by
$$\varphi_k(I+tX_1+t^2X_2+\cdots +t^lX_l) = X_k.$$
\end{cor}

\section{The Group $H^1(\slnft, V)$, $n\ge 3$}\label{h1}
Having dispensed with  the case $n=2$, we now turn our attention to
the groups $\slnft$ for $n\ge 3$.

We begin by computing the groups $H^1(\gs, V)$, where $\sigma$ is a
simplex in $\T$.  We have a split short exact sequence
$$1\lra \cs \lra \gs \stackrel{t=0}{\lra} \ps \lra 1.$$
By considering the Hochschild--Serre \sps, we see that we
have
$$H^1(\gs, V) = H^1(\ps, V) \oplus \text{Hom}_{\ps}(H_1(\cs),V).$$

In light of Proposition \ref{vanishing}, we have the following result.

\begin{prop}
If $V$ is not the adjoint representation, then $$H^1(\slnft, V) = H^1(\slnf, V).$$
\end{prop}

\begin{pf}  One checks easily in this case that $H^0(\T, {\cal H}^1) = 
H^1(\slnf, V)$
by considering the filtration $W^\bullet$ of $\T$.
\end{pf}

Now, if $V$ is the adjoint representation we have already seen that the
situation can be much different.

\begin{theorem}\label{mainthm}
If $V$ is the adjoint representation, then
$$H^1(\slnft, V) = H^1(\slnf, V) \oplus \text{\em Hom}_B(H_1(C), V).$$
\end{theorem}

\begin{pf}  As before, one checks that
$$H^0(\T, {\cal H}^1) = H^1(\slnf, V) \oplus H^0(\T, \Hom(C,V)).$$
It remains to compute the kernel of the map
$$d:C^0(\T,\Hom(C,V)) \lra C^1(\T,\Hom(C,V)).$$

Suppose that $v$ is a vertex in $\T$.  The discussion following the proof of Proposition
\ref{vanishing} shows that the group $\text{Hom}_{P_v}(H_1(C_v),V)$ has
basis the maps $\varphi_k$ where
$$\varphi_k(I+tX_1+\cdots +t^lX_l) = X_k$$
and $k\le \text{max}\{r_i-r_j\}$ (here $v=[e_1t^{r_1},e_2t^{r_2},\dots ,
e_{n-1}t^{r_{n-1}},e_n]$).  We will show that only the map $\varphi_1$
corresponds to an element in the kernel of $d$.  Note that $d$ is
$F$-linear since it is a linear combination of restriction maps, each
of which is $F$-linear.

Suppose that $k \ge 2$ and consider the map $\varphi_k$.  This map
belongs to infinitely many $\text{Hom}_{P_v}(H_1(C_v),V)$.  However,
given a vertex $v$, the map $\varphi_k$ can only belong to $\text{Hom}_{P_v}(H_1(C_v),V)$ if $r_i=r_{i+1}$ for some $i$ ({\em i.e.,}
$v$ cannot lie in $\T - W^{(n-2)}$).  In this case, there exists
a vertex $v'$ adjacent to $v$ (where $r_i > r_{i+1}$) with $\varphi_k \not\in \text{Hom}_{P_{v'}}(H_1(C_{v'}),V)$ and such that $C_e=
C_v\cap C_{v'} = C_v$.  The map $d:C^0(\T,\Hom(C,V)) \lra C^1(\T,\Hom(C,V))$ is then
$$d:(\varphi_k)_v \mapsto (\varphi_k)_e + \sum_{e'} (\varphi_k)_{e'}$$
where the sum is over the edges $e'$ incident with $v$.  However, since
$\varphi_k \not\in \text{Hom}_{P_{v'}}(H_1(C_{v'}),V)$, the $e$ component
of $d(\varphi_k)$ cannot be killed by a corresponding element of
$\text{Hom}_{P_{v'}}(H_1(C_{v'}),V)$, {\em i.e.,} $d(\varphi_k) \ne 0$.

However, the map $\varphi_1$ does belong to each $\text{Hom}_{P_v}(H_1(C_v),V)$ and to each $\text{Hom}_{P_e}(H_1(C_e),V)$.  Moreover, given an edge $e$ with
vertices $v$ and $v'$, the restriction maps $$\text{Hom}_{P_v}(H_1(C_v),V)
\lra \text{Hom}_{P_e}(H_1(C_e),V)$$
and
$$\text{Hom}_{P_{v'}}(H_1(C_{v'}),V) \lra \text{Hom}_{P_e}(H_1(C_e),V)$$
map $\varphi_1$ to the same element (with opposite sign).  It follows
that the kernel of $d:C^0 \lra C^1$ is spanned by the element $x=(x_v)_{v\in \T}$ defined by
$$x_v= \varphi_1 \in \text{Hom}_{P_v}(H_1(C_v),V).$$
If we identify this map with the map $\varphi_1 \in \text{Hom}_{B}(H_1(C),V)$, then we see that
$$H^0(\T,\Hom(C,V)) \cong \text{Hom}_B(H_1(C),V).$$
This completes the proof of Theorem \ref{mainthm}. 
\end{pf}

\section{A Second Spectral Sequence}\label{spectral2}

There is another spectral sequence which can be used to compute
the group $H^1(\slnft, V)$; namely, the Hochschild--Serre \sps 
associated to the split extension
$$1 \lra K \lra \slnft \stackrel{t=0}{\lra} \slnf \lra 1$$
where $K$ denotes the subgroup of $\slnft$ consisting of those
matrices which are congruent to the identity modulo $t$.  Since
$K$ acts trivially on $V$ and the extension splits, we see that the
map
$$d_2^{0,1}:E_2^{0,1} \lra E_2^{2,0}$$
vanishes and hence,
\begin{eqnarray*}
H^1(\slnft, V) & = & E_2^{1,0} \oplus E_2^{0,1} \\
                      & = & H^1(\slnf, V) \oplus H^0(\slnf, H^1(K,V)) \\
                      & = & H^1(\slnf, V) \oplus \text{Hom}_{\slnf}(H_1(K),V).
\end{eqnarray*}

The results of the preceding section now imply the following result.

\begin{cor}
If $V$ is not the adjoint representation, then $$\text{\em Hom}_{\slnf}
(H_1(K),V) = 0;$$ that is, there are no $\slnf$-equivariant maps
$H_1(K)\ra V$.  If $V$ is the adjoint representation, then
$$\text{\em Hom}_{\slnf}(H_1(K),V) = \text{\em
Hom}_B(H_1(C),V). \hfill \qed$$
\end{cor}

In other words, $\slnf$-equivariant maps $H_1(K) \ra V$ are in
one-to-one correspondence with $B$-equivariant maps $H_1(C) \ra V$.
Moreover, if we define $H_1^{\text{alg}}(K)$ to be the maximal algebraic
quotient of $H_1(K)$ ({\em i.e.} the product of those $\slnf$
representations which admit an equivariant homomorphism $H_1(K) \ra
V$) then we have
$$H_1^{\text{alg}}(K) = \text{Ad}.$$

The conclusion of Corollary 6.1 is clear in the case $n=2$ since we have the free product
decomposition (see {\em e.g.} \cite{knudson1})
$$K= *_{s\in {\Bbb P}^1(F)} sCs^{-1} \qquad ({\Bbb P}^1(F) = \sslf/B)$$
which implies that
$$H_1(K) = \text{Ind}_B^{\sslf} H_1(C).$$
By Shapiro's Lemma, we have 
$$H^0(B,H^1(C,V)) = H^0(\sslf, \text{Coind}_B^{\sslf} H^1(C,V));$$
that is,
\begin{eqnarray*}
\text{Hom}_B(H_1(C),V) & = & \biggl(\prod_{\sslf/B} H^1(C,V)\biggr)^{\sslf} \\
        & = &  \text{Hom}_{\sslf}(\bop_{\sslf/B}H_1(C),V) \\
        & = &  \text{Hom}_{\sslf}(H_1(K),V).
\end{eqnarray*}
Note that this argument works for any field $F$.

\begin{prop}  If $F$ is any field, then
$$H^1(SL_2(F[t],V) = H^1(SL_2(F),V) \oplus \text{\em Hom}_B(H_1(C),V).
\hfill \qed$$
\end{prop}

\section{Central Series}\label{centralseries}
In this section, we allow $F$ to be any field of any characteristic.

For each $i\ge 1$, denote by $K^i$ the subgroup of $K$ consisting of
matrices congruent to the identity modulo $t^i$.  Each $K^i$ is a normal
subgroup of $\slnft$ as it is the kernel of the map
$$\slnft \lra SL_n(F[t]/(t^i)).$$

\begin{lemma}
For each $i,j$, $[K^i,K^j] \subseteq K^{i+j}$.
\end{lemma}

\noindent {\em Proof.} If $X=I + t^iX_i + \cdots + t^lX_l$ and
$Y=I + t^jY_j + \cdots t^mY_m$, then
\begin{eqnarray*}
XYX^{-1}Y^{-1}  & = &  (I+t^iX_i +\cdots)(I+t^jY_j+\cdots)  \\
                        &  &  \hskip 1in (I-t^iX_i+\cdots)(I-t^jY_j+\cdots) \\
                    & = & (I+t^iX_i+t^jY_j +\cdots)(I-t^iX_i-t^jY_j+\cdots) \\
                   & = & I + t^{i+j}(-X_iY_j-Y_jX_i) + \cdots \\
                   & \in & K^{i+j}. \hfill \qed
\end{eqnarray*}

For each $i$ define a map $\rho_i: K^i \lra \lsln$ by
$$\rho_i(I+t^iX_i +\cdots + t^lX_l) = X_i.$$
One checks easily that $\rho_i$ is a surjective homomorphism (see
Lee and Szczarba \cite{lee} for the $\zz$ case) with kernel $K^{i+1}$.
Hence for each $i\ge 1 $, we have an isomorphism
$$K^i/K^{i+1} \cong \lsln.$$
For any group $G$, denote by $\Gamma^\bullet$ the lower central series of $G$.

\begin{lemma}
For each $i$, we have $\Gamma^i \subseteq K^i$.
\end{lemma}

\begin{pf}  Since $K^i/K^{i+1}$ is abelian for each $i$, $K^\bullet$ is
a central series in $K$ and as such contains the lower central series.
\end{pf}

Given a group $G$ and a central series $G=G^1\supseteq G^2 \supseteq 
\cdots$ we have the associated graded Lie algebra
$$\text{Gr}^\bullet G = \bop_{i\ge 1} G^i/G^{i+1}.$$

\begin{lemma} 
The  graded algebra $\text{\em Gr}^\bullet G$ is generated by
$\text{\em Gr}^1 G$ if and only if $G^l = G^{l+1}\Gamma^l$ for each
$l\ge 1$.
\end{lemma}

\begin{pf}  Note that $\text{Gr}_{\Gamma}^{\bullet} G = \bop_{i\ge 1}
\Gamma^i/\Gamma^{i+1}$ is generated by $\text{Gr}_{\Gamma}^1 G
= \Gamma^1/\Gamma^2$ and that we have a surjective map
$$\Gamma^1/\Gamma^2 \lra G^1/G^2.$$  For each $l$ we have an exact
sequence
$$0\ra (\Gamma^l\cap G^{l+1})/\Gamma^{l+1} \ra \Gamma^l/\Gamma^{l+1}
\ra G^l/G^{l+1} \ra G^l/G^{l+1}\Gamma^l \ra 0.$$
We also have a commutative diagram
$$\begin{CD}
(\Gamma^1/\Gamma^2)^{\otimes l} @>>>  (G^1/G^2)^{\otimes l}  \\
@VVV                                                        @VV{\varphi}V \\
\Gamma^l/\Gamma^{l+1}                 @>{\psi}>> G^l/G^{l+1}
\end{CD}$$

Suppose that $\text{Gr}^\bullet G$ is generated by $\text{Gr}^1 G$.  Then
$\varphi$ is surjective.  Commutativity of the diagram forces $\psi$ to
be surjective (since the top horizontal and left vertical maps are 
surjective).  But then the exact sequence
$$\Gamma^l/\Gamma^{l+1} \stackrel{\psi}{\lra} G^l/G^{l+1} 
\lra G^l/G^{l+1}\Gamma^l \lra 0$$
implies that $G^l/G^{l+1}\Gamma^l = 0$; that is, $G^l = G^{l+1}\Gamma^l$.

Conversely, if $G^l=G^{l+1}\Gamma^l$ for each $l$, then $\psi$ is
surjective.  But then $\varphi$ is also surjective; that is, $\text{Gr}^
\bullet G$ is generated by $\text{Gr}^1 G$.
\end{pf}

Now consider the algebra $\text{Gr}^\bullet K = \bop_{i\ge 1} K^i/K^{i+1}$.
For each $i$, $\text{Gr}^i K = \lsln$.  Since $\lsln = [\lsln,\lsln]$ (unless
$n=2$, $\text{char}\; F=2$), 
$\text{Gr}^\bullet K$ is generated by $\text{Gr}^1 K$.  By the lemma, we
have $K^l = K^{l+1}\Gamma^l$ for each $l$.  Note that this implies that
$K^l = K^{l+m}\Gamma^l$ for each $m\ge 1$.  Since $\bigcap_m
K^m = \{I\}$, we see that the filtrations $K^\bullet$ and
$\Gamma^\bullet$ are close together in a certain sense.  However, we
will see that in general $K^l \ne \Gamma^l$.

\medskip

\noindent {\em Remark.}  The fact that $K^l = K^{l+m}\Gamma^l$ implies that in the
completed group $\widehat{K} = \lim\limits_{\leftarrow} K/K^i$, we
have $\widehat{K}^l = \widehat{\Gamma}^l$ for each $l$.

\medskip

Consider the short exact sequence
$$0 \lra K^2/\Gamma^2 \lra H_1(K) \stackrel{{\rho_1}_*}{\lra} {\frak
sl}_n(F) \lra 0.$$
The group $K^2$ is often denoted by $SL_n(F[t],(t^2))$ and we have
$$\Gamma^2 = [SL_n(F[t],(t)),SL_n(F[t],(t))].$$

\begin{lemma}  If $F$ is a field of characteristic not equal to $2$ or $3$
or if $F$ is perfect, then $\Gamma^2 \subset E(F[t],(t^2))$.
\end{lemma}

\begin{pf}  According to Vaserstein \cite{vas}, the stable commutator
subgroup $$[GL(R,I),GL(R,J)]$$ is generated by elements
$[E_{rs}(a),E_{sr}(b)]$ with $(a,b) \in (I\times J) \cup (R\times
IJ)$ along with the $E_{rs}(x)$ with $x\in IJ$ (the latter elements
are clearly in $\Gamma^2$ and in $E(F[t],(t^2))$). The corresponding Mennicke symbol is $\left [ {a^2b \atop 1-ab}
\right ]$.  We show that these symbols vanish in $SK_1(F[t],(t^2))$.

To do this, we make use of computations in Krusemeyer
\cite[Lemma 12.3]{kruse}.  The hypotheses on $F$ are needed to show the
relations $\left [ {t^2 \atop 1-t^2k} \right ] = 1$ and $\left [ {t^3
\atop 1-t^2k} \right ] = 1$ for $k \in F[t]$.  The hardest case is
where $(a,b)\in (t)\times (t)$.  Write $a=tf, b=tg, g=g_0 + g_1t +
t^2h$, where $g_i \in F$ and $f,h\in F[t]$.  Then
\begin{eqnarray*}
\left [ {a^2b \atop 1-ab} \right ] &  =  &  \left [ {t^3f^2g \atop
1-t^2fg} \right ] = \left [ {t^3f^2g \atop 1-t^2fg} \right ] \left [
{t^2 \atop 1-t^2fg} \right ] \\
  &  =  & \left [ {t^5f^2g \atop 1-t^2fg} \right ] = \left [ {t^3f
\atop 1-t^2fg} \right ] \left [ {fg \atop 1-t^2fg} \right ] \\
  &  =  & \left [ {t^3f \atop 1-t^2fg} \right ] = \left [ {t^5f \atop
  1-t^2fg} \right ] = \left [ {t^2f \atop 1-t^2fg} \right ] \\
  &  =  & \left [ {t^2f \atop 1-t^3fg_1} \right ] = \left [ {t^3f
  \atop 1-t^3fg_1} \right ] = 1.
\end{eqnarray*}
Similar computations show that $\left [ {a^2b \atop 1-ab} \right ] =
1$ for $(a,b) \in F[t]\times (t^2)$.  This completes the proof.
\end{pf}

The lemma implies that we have a surjective map
$$K^2/\Gamma^2 \lra SK_1(F[t],(t^2)).$$  The latter group was computed
by Krusemeyer \cite{kruse}; it equals the module of differentials
$\Omega_F^1$, if $\text{char}\, F \ne 2,3$ or if $F$ is perfect.
Observe that this is nonzero in general.

\begin{prop}  Suppose $\text{char}\, F \ne 2,3$ or $F$ is perfect.  If $\Omega_F^1 \ne 0$, then $K^l \ne \Gamma^l$ for all
$l\ge 2$.
\end{prop}

\begin{pf}  Since $K^2/\Gamma^2$ surjects onto $\Omega_F^1$, we have
$K^2 \ne \Gamma^2$.  Suppose $K^l = \Gamma^l$ for some $l$.  Then
$K^2 = K^l\Gamma^2 = \Gamma^l\Gamma^2 \subseteq \Gamma^2$, a
contradiction.
\end{pf}

If $\Omega_F^1 = 0$, then we get no information about the structure of
$K^2/\Gamma^2$.  Note that if $F$ is a finite field, then $\Omega_F^1
= 0$.

\section{The Action of $K$ on $\X$}\label{action}

We now turn our attention to the computation of the group $H_1(K)$ in
the case where $F$ is a finite field.  While much of what we say below
is true for any field, we restrict our attention to the finite case.

In \cite[Prop. 4.1]{knudson2}, we found a fundamental domain for the action of $K$
on $\X$.  Denote by $S$ a set of coset representatives for $\slnf/B$,
where $B$ is the upper triangular subgroup.  Define a subcomplex $\Z$ of
$\X$ by
$$\Z = \bigcup_{s\in S} s\T.$$
Then $\Z$ is a fundamental domain for the action of $K$ on $\X$.

If $\sigma$ is a simplex in $\Z$, denote by $\ks$ the stabilizer of
$\sigma$ in $K$.  Since $\X$ is contractible, we have a spectral sequence
converging to $H_\bullet(K)$ with $E^1$--term
\begin{equation}\label{ssk}
E^1_{p,q} = \bop_{\dim \sigma = p} H_q(\ks)
\end{equation}
where $\sigma$ ranges over the simplices of $\Z$.

Observe that the bottom row $E^1_{*,0}$ is simply the simplicial chain
complex $S_\bullet(\Z)$.  Note that $\Z$ is contractible since the
index set $S$ is finite ($F$ is a finite field) and each $s\in S$
fixes the initial vertex $v_0$.  Moreover, $s\T \cap s'\T$ is at most
a codimension one face.  It follows that the straight line contracting
homotopies for each $s\T$ can be glued together to obtain a
contracting homotopy of $\Z$.  It follows that
$$E^2_{p,0} = \begin{cases}
                              \zz    &   p=0 \\
                                 0    &  p>0.
                      \end{cases}$$
Thus, to compute $H_1(K)$ we need only compute the group $E^2_{0,1}$.

\section{Reduction}\label{reduction}

For each $q$ define a coefficient system $\cH_q$ on $\Z$ by setting
$\cH_q(\sigma) = H_q(\ks)$.  Then the $q$th row of the spectral sequence
(\ref{ssk}) is simply the chain complex $C_\bullet(\Z,\cH_q)$, and
the group $E^2_{p,q}$ is the group $H_p(\Z,\cH_q)$.

In the previous section, we showed that $H_1(K) = E^2_{0,1}$.  Thus we
have $$H_1(K) = H_0(\Z,\cH_1).$$  To compute the latter group, we first
define a sequence of subcomplexes of $\Z$.

Denote by $\Z^{(1)}$ the 1-skeleton of $\Z$.  Clearly, $H_0(\Z,\cH_1)
= H_0(\Z^{(1)},\cH_1)$.  For each $i\ge 1$, let $\Z_i$ be the subgraph of
$\Z^{(1)}$ spanned by the vertices having distance at most $i$ from
$v_0$.  Then  $$\Z^{(1)} = \bigcup_{i\ge 1}\Z_i$$
and
$$H_0(\Z^{(1)},\cH_1) = \varinjlim H_0(\Z_i,\cH_1).$$

\begin{prop}
The inclusion $\Z_1 \lra \Z^{(1)}$ induces a surjection
$$H_0(\Z_1,\cH_1) \lra H_0(\Z^{(1)},\cH_1).$$
\end{prop}

\begin{pf}  Consider the long exact sequence of the pair $(\Z_{i+1},\Z_i)$
for $i\ge 1$:
$$H_1(\Z_{i+1},\Z_i;\cH_1) \ra H_0(\Z_i,\cH_1) \ra H_0(\Z_{i+1},\cH_1) \ra H_0(\Z_{i+1},\Z_i;\cH_1) \ra 0.$$
Let $v$ be a vertex in $\Z_{i+1}$.  Then there exists an $s \in \slnf/B$ and
a vertex $v' \in \T$ such that $v=sv'$.  It follows that $K_v=sK_{v'}s^{-1}$
 and that $H_1(\kv)$ admits a 
decomposition
$$H_1(\kv) = tW_1 \oplus t^2W_2 \oplus \cdots \oplus t^{i+1}W_{i+1}$$
where each $W_k$ is an $F$--vector space spanned by certain elements of
$\lsln$ (this follows from the decomposition of $H_1(K_{v'}) = H_1(C_{v'})$ described in Section \ref{psigma}).

Let $v\in \Z_{i+1}-\Z_i$.  Then there exist $v_1,\dots ,v_l \in \Z_i$
adjacent to $v$ with $K_{v_j} \subset \kv$ for each $j$.  If $e_j$ denotes
the edge connecting $v_j$ with $v$, then $K_{e_j} = K_{v_j}$.  Now, the
boundary map $$\partial:C_1(\Z_{i+1},\Z_i;\cH_1) \ra C_0(\Z_{i+1},\Z_i;
\cH_1)$$ maps $\displaystyle\bop_{j=1}^l H_1(K_{e_j})$ surjectively 
onto the summand $tW_1 \oplus \cdots \oplus t^iW_i$ of
$$H_1(\kv) = tW_1 \oplus \cdots \oplus t^iW_i \oplus t^{i+1}W_{i+1}.$$
Figure 2 shows an example consisting of the vertices
$v=[t^2e_1,e_2,t^2e_3]$, $v'=[t^2e_1,e_2,te_3]$,
$v''=[te_1,e_2,t^2e_3]$, $v_1=[te_1,e_2,te_3]$, where the top group is
$H_1(C_v)$, the three middle groups are $H_1(C_{v'',v})$,
$H_1(C_{v_1,v})$, $H_1(C_{v',v})$, and the bottom groups are
$H_1(C_{v''})$, $H_1(C_{v_1})$, $H_1(C_{v'})$.  The space $W_{ij}$
is a one-dimensional weight space for the root $\alpha_{ij}$.
The vertices $v$,
$v'$, $v''$ lie in $\Z_2$ and $v_1$ lies in $\Z_1$.

\begin{figure}
$$\begin{array}{ccccccccc}
            &         &           &          & tW_{12}\oplus &  &  &
            &   \\
            &         &           &          &t^2W_{12}\oplus & & & 
            &   \\
            &         &           &          & tW_{32}\oplus & & & 
            &   \\
            &         &           &          &t^2W_{32} & & & 
            &   \\
            &         &           &\nearrow  & \uparrow &\nwarrow & &
            &   \\
            &         &tW_{12}\oplus &       &tW_{12}   & &tW_{12}\oplus & 
            &   \\
            &         &tW_{32}\oplus &       &\oplus    & &t^2W_{12}\oplus &
            &   \\
            &         &t^2W_{32}     &       &tW_{32}   & &tW_{32} & 
            &   \\
            &\swarrow &              &       &\downarrow & &       &
\searrow    &   \\
tW_{12}\oplus &       &              &       &tW_{12}   & &        &
            &tW_{12}\oplus \\
tW_{31}\oplus &       &              &       &\oplus    & &        &
            &tW_{13}\oplus \\
tW_{32}     &         &              &       &tW_{32}   & &        &
            &tW_{32}       
\end{array}$$
\caption{The Vanishing of $H_0(\Z_{i+1},\Z_i;\cH_1)$}
\end{figure}

It follows that the group $H_0(\Z_{i+1},\Z_i;\cH_1)$ can consist only of
classes arising from the various $t^{i+1}W_{i+1}$.  Moreover, such a
class $t^{i+1}w_{i+1}$ can arise only from vertices of the form $v=s\dot [e_1t^{r_1},
\dots , e_{n-1}t^{r_{n-1}},e_n]= sv'$ where $r_k=r_{k+1}$ for some $k=
1,2,\dots n-1$ (set $r_n=0$).  (Recall that those vertices for which
the $r_k$ are positive and distinct satisfy $H_1=tW_1\oplus \cdots 
\oplus t^iW_i$---see the discussion following Proposition 4.2.)
  However, $v$ is adjacent to a vertex $w$ in $\Z_{i+1}
- \Z_i$ for which $r_k>r_{k+1}$ and hence the corresponding class
$t^{i+1}w_{i+1}$ is trivial in $H_1(K_w)$.  Also, the class $t^{i+1}w_{i+1}$
is nontrivial in $H_1(K_e)$ where $e$ is the edge joining $v$ and $w$.
It follows that the boundary map hits the class $t^{i+1}w_{i+1}$ in
$H_1(\kv)$; that is, $H_0(\Z_{i+1},\Z_i;\cH_1) = 0$.
Figure 2 illustrates this for the classes $t^2w_{12}, t^2w_{32} \in
H_1(C_v)$.  The class $t^2w_{12}$ is hit by the class $t^2w_{12} \in
H_1(C_{v'',v})$ (notice that this class is trivial in $H_1(C_{v''})$).
Similarly, the class $t^2w_{32}$ is hit by $t^2w_{32}\in H_1(C_{v',v})$. 
Since $H_0(\Z^{(1)},\cH_1)=\varinjlim H_0(\Z_i,\cH_1)$ and each map
$H_0(\Z_i,\cH_1)\ra H_0(\Z_{i+1},\cH_1)$ is surjective, we have a
surjection $H_0(\Z_1,\cH_1)\ra H_0(\Z^{(1)},\cH_1)$.
\end{pf}

\medskip

\noindent {\em Remark.}  It is crucial that $n\ge 3$ in the above proof.
When $n=2$, each vertex $v_{i+1}$ in $\Z_{i+1}-\Z_i$ is adjacent to a
single vertex $v_i$ in $\Z_i$. The stabilizer $C_{e_i}$ of the edge
joining $v_i$ to $v_{i+1}$ maps 
surjectively onto the part of $C_{v_{i+1}}$ consisting of polynomials of
degree less than $i+1$.  However, $v_{i+1}$ is not adjacent to {\em any}
vertex in $\Z_{i+1}-\Z_i$.  Hence, $$H_0(\Z_{i+1},\Z_i;\cH_1)=
\bop_{s\in \sslf/B} F$$ with basis the various $t^{i+1}w_{i+1}$.

\medskip

\begin{cor}  The group $H_1(K)$ is a quotient of $H_0(\Z_1,\cH_1)$.
\end{cor} 

If the field $F$ is finite, then we see immediately that the graph $\Z_1$
is finite and hence $H_1(K)$ is finite dimensional.  This is in sharp
contrast to the case $n=2$. 

\section{Computation}\label{computation}

We are now in a position to compute the group $H_1(K)$ in certain cases.
 
\begin{theorem}  If $n=3$ and $F={\Bbb F}_2$, then $H_1(K) = {\frak
sl}_3({\Bbb F}_2)$.
\end{theorem}

\begin{pf}  By Proposition 9.1, we have a surjective map
$$H_0(\Z_1,{\cal H}_1) \lra H_1(K).$$
Since $H_1(K)$ surjects onto ${\frak sl}_3(F)$, we need only check that
$H_0(\Z_1,{\cal H}_1)$ is an 8-dimensional vector space.

The
complex $\Z_1$ is the incidence geometry of ${\Bbb F}_2^3$.  It is a
graph with 14 vertices and 21 edges.  For each vertex $v$ in $\Z_1$,
the group $H_1(\kv)$ is two-dimensional and for each edge $e$, the
group $H_1(K_e)$ is one-dimensional.  It follows that
$$\partial:C_1(\Z_1,{\cal H}_1) \lra C_0(\Z_1,{\cal H}_1)$$
is a map from a 21-dimensional vector space to a 28-dimensional
vector space. Moreover, this map is clearly ${\Bbb F}_2$-linear.     
Denote the generators of $C_0(\Z_1,{\cal H}_1)$
by $f_1,f_2,\dots ,f_{28}$.  The relations imposed by $\partial$ are
as follows:
$$\begin{array}{lll}
f_1=f_3         &  f_{15}=f_2+f_4         &  f_{14}=f_{20}           \\
f_2=f_{11}     & f_{17}=f_5+f_6          & f_{27}=f_{13}+f_{14} \\
f_{12}=f_{10} & f_{19}=f_7+f_8          & f_{27}+f_{28}=f_{17}+f_{18}\\
f_9=f_7          & f_{21}=f_9+f_{10}      & f_{28}=f_{21}+f_{22}  \\
f_5=f_8          & f_{23}=f_{11}+f_{12}  & f_{25}+f_{26}=f_{19}+f_{20}\\
f_4=f_6          & f_{18}=f_{24}             & f_{26}=f_{23}+f_{24}  \\
f_{13}=f_1+f_2 & f_{16}=f_{22}          & f_{25}=f_{15}+f_{16}
\end{array}$$
By considering the $28 \times 21$ matrix associated to $\partial$, we
see that $H_0(\Z_1,{\cal H}_1)$ is an 8-dimensional $F$-vector space
(the rank of $\partial$ is 20).
\end{pf}

\noindent {\em Remarks.} 1.  A similar argument can be used to compute
$H_1(K) = {\frak sl}_3({\Bbb F}_3)$ for $F={\Bbb F}_3$.  In this case
the complex $\Z_1$ has 25 vertices and 42 edges.  The map $\partial$
is injective and hence $H_0(\Z_1,\cH_1)$ is 8-dimensional.

\medskip

\noindent 2.  To prove that $H_1(K)=\lsln$ for all $n$,
it suffices to show that the vector space  $H_0(\Z_1,{\cal H}_1)$     
has dimension $n^2-1$.
Since $\Z_1$ is a finite complex it should be possible to carry this out.
However, the combinatorics seem to be rather complicated in general.
The complexity increases as $n$ does and for a given $n$, the
computation becomes more difficult as the cardinality of $F$ increases.

\medskip

\noindent 3.  As a corollary, we see that the subgroup
$K^2$ of $K$ consisting of matrices congruent to the identity modulo $t^2$
 is equal to the commutator subgroup of $K$ for $n=3$, $F={\Bbb F}_2,
{\Bbb F}_3$.


\begin{thebibliography}{99}

\bibitem{tits} F. Bruhat, J. Tits,  Groupes r\'eductifs sur un corps local
I: Donn\'ees radicielles valu\'ees, Publ. IHES {\bf 41} (1972), 5--252.

\bibitem{dwyer} W. Dwyer, Twisted homological stability for general
linear groups, Ann. of Math. {\bf 111} (1980), 239--251.

\bibitem{fulton} W. Fulton, J. Harris, ``Representation Theory:  A First 
Course'', Springer--Verlag, Berlin/Heidelberg/New York, 1991.

\bibitem{humphreys} J. Humphreys, ``Linear Algebraic Groups'', Springer--Verlag, Berlin/Heidelberg/New York, 1975.

\bibitem{vdk} W. van der Kallen, Homology stability for linear groups,
Invent. Math. {\bf 60} (1980), 269--295.

\bibitem{knudson1} K. Knudson,  The homology of $SL_2(F[t,t^{-1}])$, J.  Alg. {\bf 180} (1996), 87--101.

\bibitem{knudson2} K. Knudson, The homology of special linear
groups over polynomial rings, Ann. Sci. \'Ecole Norm. Sup. (4) {\bf
30}, no. 3 (1997), 385--415.

\bibitem{kruse} M. Krusemeyer, Fundamental groups, algebraic
$K$-theory, and a problem of Abhyankar, Invent. Math. {\bf 19} (1973), 15--47. 

\bibitem{lee} R. Lee, R. Szcarba,  On the homology and cohomology of
congruence subgroups, Invent. Math. {\bf 33} (1976), 15-53.

\bibitem{nagao} H. Nagao, On $GL(2,K[x])$, J. Poly. Osaka Univ.
{\bf 10} (1959), 117--121.

\bibitem{serre} J.-P. Serre,``Trees'', Springer--Verlag,
Berlin/Heidelberg/New York, 1980.

\bibitem{soule1} C. Soul\'e, Chevalley groups over polynomial rings, in Homological Group Theory (C.T.C. Wall, ed.),  London Math. Soc. Lecture Notes {\bf 36}, Cambridge University Press, Cambridge (1979), 359--367.

\bibitem{soule2} C. Soul\'e, The cohomology of $SL_3(\zz )$, Topology {\bf 17} (1978),
1--22.

\bibitem{vas} L. Vaserstein, Foundations of algebraic $K$-theory,
Russian Math. Surveys {\bf 31} (1976), 89--156.

\end{thebibliography}
\end{document}